\documentclass{amsart}
\usepackage{epsfig}
\theoremstyle{plain}
\newtheorem{theorem}{Theorem}[section]
\newtheorem{lemma}[theorem]{Lemma}
\newtheorem{proposition}[theorem]{Proposition}
\newtheorem{corollary}[theorem]{Corollary}

\theoremstyle{definition}

\newtheorem{remark}[theorem]{Remark}

\theoremstyle{remark}
\newtheorem*{remark*}{Remark}

\def\bd{\partial}
\def\C{{\mathbb C}}
\def\R{{\mathbb R}}

\def\H{{\mathbb H}}
\def\Z{{\mathbb Z}}
\def\g{{\mathfrak g}}
\def\lap{\triangle}

\def\arctanh{{\mathop{\rm arctanh}}}
\def\area{\mathop{\rm area}}
\def\grad{{\nabla}}

\def\M{{X}} 
\def\N{{M}} 
\def\b{{b}} 
\def\l{{\ell}} 

\def\cal{\mathcal}

\def\DS{{\cal H} {\cal D} {\cal S}}

\title{Harmonic deformations of hyperbolic 3-manifolds}
\author{Craig D. Hodgson}
\address{Department of Mathematics and Statistics\\
University of Melbourne\\
Victoria 3010, Australia}
\email{cdh@ms.unimelb.edu.au}
\thanks{The research of the first author was partially
supported by grants from the ARC}
\author{Steven P. Kerckhoff}
\address{Department of Mathematics\\
Stanford University\\
Stanford, CA 94305, U.S.A.}
\email{spk@math.stanford.edu}
\thanks{The research of the second author was partially
supported by grants from the NSF}

\begin{document}

\begin{abstract}
This paper gives an exposition of the authors' harmonic deformation
theory for 3-dimensional hyperbolic cone-manifolds. We discuss
topological applications to hyperbolic Dehn surgery as well
as recent applications to Kleinian group theory. A central idea
is that local rigidity results (for deformations fixing cone angles)
can be turned into effective control on the deformations that do exist.
This leads to precise analytic and geometric versions of the idea
that hyperbolic structures with short geodesics are close to
hyperbolic structures with cusps. The paper also outlines
a new harmonic deformation theory which applies whenever there is
a sufficiently large embedded tube around the singular locus,
removing the previous restriction to cone angles at most $2\pi$.
\end{abstract}

\maketitle

\section{Introduction}\label{intro}

The local rigidity theorem of Weil \cite{weil} and Garland \cite{Ga} 
for complete,
finite volume hyperbolic
manifolds states that there is no non-trivial deformation of
such a structure through {\it complete} hyperbolic structures
if the manifold has dimension at least 3.  If the manifold is
closed, the condition that the structures be complete is automatically
satisfied.  However, if the manifold is non-compact, there may be
deformations through incomplete structures.  This cannot
happen in dimensions greater than 3 (Garland-Raghunathan \cite{GR}); but
there are always non-trivial deformations in dimension 3 (Thurston 
\cite{thnotes})
in the non-compact case.

In \cite{HK1} this rigidity theory is extended to
a class of finite volume, orientable 3-dimensional hyperbolic
cone-manifolds, {\it i.e.} hyperbolic structures on 3-manifolds with
cone-like singularities along a knot or link.
The main result is that such structures are locally rigid if the
cone angles are fixed, under the extra hypothesis that
all cone angles are at most $2\pi$.
There is a smooth, incomplete structure on the complement of the
singular locus; by completing the metric the singular cone-metric is recovered.
The space of deformations of (generally incomplete) hyperbolic structures
on this open manifold has non-zero dimension, so
there will be deformations if the cone angles are allowed to vary.
An application of the implicit function theorem shows that it is possible
to deform the structure so that the
metric completion is still a 3-dimensional hyperbolic cone-manifold, 
and it is always
possible to deform the cone-manifold to make arbitrary (small) changes
in the cone angles.  In fact, the collection of cone angles locally 
parametrizes
the set of cone-manifold structures.

A (smooth) finite volume hyperbolic $3$-manifold with cusps
is the interior of a compact $3$-manifold with torus boundary 
components.  Filling
these in by attaching solid tori produces a closed manifold; there is 
an infinite
number of topologically distinct ways to do this, parametrized by the isotopy
classes of the curves on the boundary tori that bound disks in the solid tori.
These curves are called the ``surgery curves".
The manifold with cusps can be viewed as a cone-manifold structure 
with cone angles $0$
on any of these closed manifolds.  If it is possible to increase the cone angle
from $0$ to $2 \pi$, this constructs a smooth hyperbolic structure on 
this closed
manifold.  This process is called {\em hyperbolic Dehn surgery}.
Thurston (\cite{thnotes}) proved that hyperbolic Dehn surgery fails
for at most a finite number of choices of surgery curves on each 
boundary component.

The proof of local rigidity puts strong constraints on those deformations
of hyperbolic cone-manifolds that {\em do} exist.  It is possible
to control the change in the geometric structure when the cone angles are
deformed a fixed amount.  Importantly, this control depends only on the
geometry in a tubular neighborhood around the singular locus, not on the rest
of the $3$-manifold.  In particular, it provides geometric and analytic
control on the hyperbolic Dehn surgery process.   This idea is developed
in \cite{HK2}.

That paper provides a quantitative version of Thurston's hyperbolic Dehn
surgery theorem. Applications include the first universal bounds on the
number of non-hyperbolic Dehn fillings on a cusped hyperbolic 3-manifold,
and estimates on the changes in volume and core geodesic length during
hyperbolic Dehn filling.

The local rigidity theory of \cite{HK1} was generalized by Bromberg 
(\cite{Br1})
to include geometrically finite hyperbolic cone-manifolds.  Recently,
there have been some very imaginative and interesting applications of
the deformation theory of geometrically finite hyperbolic cone-manifolds
to well-known problems in Kleinian groups.
In particular, the reader is referred to \cite{BB2} in this volume
for a description of some of these results and references to others.

Our purpose here is to provide a brief outline of the main ideas and
results from \cite{HK1} and \cite{HK2} and how they
are related to the Kleinian group applications.  As noted above,
the central idea is that rigidity results can be turned into effective
control on the deformations that do exist.  However, we wish to emphasize a
particular consequence that provides the common theme between 
\cite{HK2} and the
Kleinian group applications in \cite{Br2}, \cite{Br3}, \cite{BB1},
\cite{BB2} and \cite{BBES}.  As a corollary of the control provided by
effective rigidity, it is possible to give precise analytic and 
geometric meaning
to the familiar idea that hyperbolic structures with short geodesics 
are ``close"
to ones with cusps.   Specifically, it can be shown that a structure
with a sufficiently short geodesic can be deformed through hyperbolic
cone-manifolds to a complete structure, viewed as having cone angle $0$.
Furthermore, the total change in the structure can be proved to be
arbitrarily small for structures with arbitrarily short geodesics.
Most importantly this control is independent of the manifolds involved,
depending only on the lengths and cone angles.
 
There are varied reasons for wanting to find such a family of cone-manifolds.
It is conjectured that any closed hyperbolic $3$-manifold can be obtained
by hyperbolic Dehn surgery on some singly cusped, finite volume hyperbolic $3$-manifold.
If this were true, it could have useful implications in 
$3$-dimensional topology.
In \cite{HK2} it is proved that it is true for any such closed $3$-manifold
whose shortest geodesic has length at most $0.162$.  (See Theorem \ref{connect}
in Section \ref{dehn}.)

In \cite{Br2} Bromberg describes a construction that, remarkably, allows one
to replace an incompressible, geometrically infinite end with a short
geodesic by a geometrically finite one by gluing in a wedge that creates
a cone angle of $4 \pi$ along the short geodesic.  Pushing the cone angle
back to $2 \pi$ provides an approximation of the structure with a
geometrically infinite end by one with a geometrically finite end.  In
\cite{Br2} and \cite{BB1}, Bromberg and then Brock and Bromberg give
proofs of important cases of the Density Conjecture in this way.
This work is described in \cite{BB2}.

In general, a sequence of Kleinian groups with geodesics that are becoming
arbitrarily short (or a single Kleinian group with a sequence of arbitrarily
short geodesics) is very difficult to analyze.  Things are often simpler
when the lengths are actually {\em equal} to $0$; i.e., when they are
cusps.  Thus, if structures with short geodesics can be uniformly compared
with structures where they have become cusps, this can be quite useful.
One example where this idea has been successfully employed is \cite{BBES}.
It seems likely that there will be others in the near future.

Note that the application to the Density Conjecture above involves
cone angles between $4 \pi$ and $2 \pi$ whereas the theory in \cite{HK1}
and \cite{Br1} requires cone angles to be at most $2 \pi$.  Thus,
this application actually depends on a new version of the deformation
theory (\cite{HK3}) which applies to all cone angles, as long as there
is a tube of a certain radius around the singular locus.

Because of its connection with these Kleinian group results, we use the
current paper as an opportunity to outline the main points in this new theory.
It is based on a boundary value problem which is used to
construct infinitesimal deformations with the same essential properties
as those utilized in \cite{HK1} and \cite{HK2}.  Explaining those
properties and how they are used occupies Sections \ref{harmdef} to \ref{dehn}.
The discussion of Kleinian groups and the new deformation theory are both
contained in Section \ref{boundaryvalues}, the final section of the paper.

The relaxation of the cone angle restriction has implications for
hyperbolic Dehn surgery.  Some of these are also described in the
final section.  (See Theorems \ref{bdrylocalparam} and \ref{shapethm}.)

\section{Deformations of hyperbolic structures}\label{harmdef}

A standard method for analyzing families of structures or maps is
to look at the infinitesimal theory where the determining equations
simplify considerably.  To this end, we first describe precisely what
we mean by a $1$-parameter family of hyperbolic structures on a manifold.
Associated to the derivative of such a family are various analytic, algebraic,
and geometric objects which play a central role in this theory.  It is useful
to be able to move freely among the interpretations provided by these objects
and we attempt to explain the relationships between them.

The initial portion of this analysis is quite general, applying to hyperbolic
structures in any dimension or, even more broadly, to structures modeled on
a Lie group acting transitively and analytically on a manifold.

A hyperbolic structure on an $n$-manifold $\M$
is determined by local charts modeled on $\H^n$ whose overlap
maps are restrictions of global isometries of $\H^n$.  These
determine, via analytic continuation,
a map $\Phi :\tilde \M \to \H^n$ from the universal cover
$\tilde \M$ of $\M$ to $\H^n$, called the {\em developing map},
which is determined uniquely up to post-multiplication by an
element of $G = {\rm isom} (\H^n)$. The developing map satisfies
the equivariance property $\Phi(\gamma m) = \rho(\gamma)\Phi(m)$,
for all $m\in \tilde \M$, $\gamma \in \pi_1(\M)$,
where $\pi_1(\M)$ acts on $\tilde \M$ by covering transformations,
and $\rho: \pi_1(\M) \to G$
is the {\em holonomy representation} of the structure.
The developing map also determines the hyperbolic metric on $\tilde \M$
by pulling back the hyperbolic metric on $\H^n$.
(See \cite{thbook} and \cite{Rat} for a complete discussion of these ideas.)

We say that two hyperbolic structures are {\em equivalent}
if there is a diffeomorphism $f$ from $\M$ to itself taking one
structure to the other.  We will use the term ``hyperbolic
structure" to mean such an equivalence class.
A {\em 1-parameter family}, $\M_t$, of hyperbolic structures defines
a 1-parameter family of developing maps
$\Phi_t:\tilde \M \to \H^n$,
where two families are equivalent under the relation
$\Phi_t \equiv k_t \Phi_t \tilde f_t$
where $k_t$ are isometries of $\H^n$ and $\tilde f_t$ are
lifts of diffeomorphisms $f_t$ from $\M$ to itself.
We assume that $k_0$ and $\tilde f_0$ are the identity,
and denote $\Phi_0$ as $\Phi$.
All of the maps here are assumed to be smooth and to vary
smoothly with respect to $t$.

The tangent vector to a smooth family of hyperbolic
structures will be called an {\em infinitesimal deformation}.
The derivative at $t = 0$ of a 1-parameter family of developing
maps $\Phi_t:\tilde \M \to \H^n$ defines a map
$\dot \Phi: ~ \tilde \M \to T\H^n$.  For any point
$m \in \tilde \M$, $\Phi_t(m)$ is a curve
in $\H^n$ describing how the image of $m$ is moving under the
developing maps; $\dot \Phi (m)$ is the initial tangent vector to the
curve.

We will identify $\tilde \M$ locally with $\H^n$ and $T \tilde \M$
locally with $T\H^n$ via the initial developing map $\Phi$.
Note that this identification is generally not a global diffeomorphism unless
the hyperbolic structure is complete.  However, it is a {\em local}
diffeomorphism, providing identification of small open sets
in $\tilde \M$ with ones in $\H^n$.

In particular, each point $m \in \tilde \M$ has a neighborhood $U$ where
$\Psi_t = \Phi^{-1} \circ \Phi_t : U \to \tilde \M$ is defined,
and the derivative at $t=0$ defines a vector field on $\tilde \M$,
$v = \dot \Psi : \tilde \M \to T \tilde \M$.
This vector field determines the infinitesimal variation in developing maps
since $\dot \Phi = d\Phi \circ v$, and also determines
the infinitesimal variation in metric as follows. Let $g_t$ be the
hyperbolic metric on $\tilde \M$ obtained by pulling back the
hyperbolic metric on $\H^n$ via $\Phi_t$ and put $g_0=g$. Then
$g_t = \Psi_t^* g$ and the infinitesimal variation in metrics
$\dot g = {d  g_t \over dt} |_{t=0}$ is
the Lie derivative, ${\cal L}_v g$, of the initial metric $g$ along $v$.

Riemannian covariant differentiation of
the vector field $v$ gives a $T \tilde \M$ valued
1-form on $\tilde \M$, $\nabla v : T\tilde \M \to T \tilde \M$,
defined by $\nabla v(x) = \nabla_x v$ for $x \in T\tilde \M$.
We can decompose $\nabla v$ at each point into a symmetric part
and skew-symmetric part.
The {\em symmetric part},
$\tilde\eta = (\nabla v)_{sym}$,
represents the infinitesimal change in metric, since
$$\dot g (x,y) = {\cal L}_v g (x,y) =
g(\nabla_{\!x} v,y) + g(x,\nabla_{\!y} v) = 2 g(\tilde \eta(x),y)$$
for $x,y \in T \tilde \M$. In particular, $\tilde\eta$
descends to a well-defined $T\M$-valued 1-form $\eta$ on $\M$.
The {\em skew-symmetric} part $(\nabla v)_{skew}$ is the
{\em curl} of the vector field $v$; its value at $m \in \tilde \M$
describes the infinitesimal rotation about $m$ induced by $v$.

To connect this discussion of infinitesimal deformations
with cohomology theory, we consider the
Lie algebra  $\g$ of $G = {\rm isom} (\H^n)$
as vector fields on $\H^n$ representing infinitesimal
isometries of $\H^n$. Pulling back these vector fields
via the initial developing map $\Phi$
gives locally defined infinitesimal
isometries on $\tilde \M$ and on $\M$.

Let $\tilde E, E$ denote the vector bundles over
$\tilde \M, \M$ respectively of (germs of)
infinitesimal isometries. Then we can regard $\tilde E$ as the product
bundle with total space $\tilde \M \times {\g}$,
and $E$ is isomorphic to $(\tilde \M \times {\g})/\!\!\sim$
where $(m,\zeta) \sim (\gamma  m, Ad \rho(\gamma) \cdot \zeta)$
with $\gamma \in \pi_1(\M)$ acting on $\tilde \M$
by covering transformations and on ${\g}$
by the adjoint action of the holonomy $\rho(\gamma)$.
At each point $p$ of $\tilde \M$, the fiber of
$\tilde E$ splits as a direct
sum of infinitesimal pure translations and infinitesimal pure
rotations about $p$; these can be identified with
$T_p \tilde\M$ and $so(n)$ respectively.
The hyperbolic metric on $\tilde\M$ induces a metric on 
$T_p \tilde\M$ and on $so(n)$.  A metric can then be 
defined on the fibers of $\tilde E$ in which the two factors
are orthogonal;  this descends to a metric on the fibers of $E$.

Given a vector field $v: \tilde \M \to T \tilde \M$, we
can lift it to a section $s : \tilde \M \to \tilde E$
by choosing an ``osculating'' infinitesimal isometry $s(m)$ which best
approximates the vector field $v$ at each point $m \in \tilde \M$.
Thus $s(m)$ is the unique infinitesimal isometry
whose translational part and rotational part
at $m$ agree with the values of $v$ and $curl~v$ at $m$.
(This is the ``canonical lift'' as defined in \cite{HK1}.)
In particular, if $v$ is itself an infinitesimal isometry of
$\tilde \M$ then $s$ will be a constant section.

Using the equivariance property of the developing maps
it follows that $s$ satisfies an ``automorphic'' property:
for any fixed $\gamma \in \pi_1(\M)$, the difference
$s(\gamma m) - Ad \rho(\gamma) s(m)$ is a {\em constant} 
infinitesimal isometry,
given by the variation $\dot \rho(\gamma)$ of holonomy isometries
$\rho_t(\gamma) \in G$ (see Prop 2.3(a) of \cite{HK1}).
Here $\dot \rho: \pi_1(\M) \to \g$
satisfies the cocyle condition $\dot \rho(\gamma_1 \gamma_2) = \dot
\rho(\gamma_1)+ Ad \rho (\gamma_1) \dot \rho(\gamma_2)$,
so it represents a class in group cohomology
$[\dot \rho] \in H^1(\pi_1(\M);Ad \rho)$,
describing the variation of holonomy representations $\rho_t$.

Regarding  $s$ as a vector-valued function with values
in the vector space $\g$, its differential $\tilde\omega = ds$
satisfies $\tilde\omega(\gamma m) =  Ad \rho(\gamma) \tilde\omega(m)$ so
it descends to a closed 1-form $\omega$ on $\M$ with values in the bundle $E$.
Hence it determines a de Rham cohomology
class $[\omega]\in H^1(\M; E)$.
This agrees  with the group cohomology class $[\dot \rho]$
under the de Rham isomorphism $H^1(\M; E) \cong H^1(\pi_1(\M);Ad \rho)$.
Also, we note that the translational part of $\omega$
can be regarded as
a $T\M$-valued $1$-form on $\M$.  Its symmetric part is exactly the
form $\eta$ defined above (see Prop 2.3(b) of \cite{HK1}),
describing the infinitesimal change in metric on $\M$.

On the other hand, a family of hyperbolic structures determines
only an equivalence class of families of developing maps and
we need to see how replacing one family by an equivalent family
changes both the group cocycle and the de Rham cocycle.  Recall
that a family equivalent to $\Phi_t$ is
of the form $k_t \Phi_t \tilde f_t$ where $k_t$ are isometries of $\H^n$ and
$\tilde f_t$ are lifts of diffeomorphisms $f_t$ from $\M$ to itself.
We assume that $k_0$ and $\tilde f_0$ are the identity.

The $k_t$ term changes the path $\rho_t$ of holonomy representations
by conjugating by $k_t$.  Infinitesimally, this changes the
cocycle $\dot \rho$ by a coboundary in the sense of group cohomology.
Thus it leaves the class in $H^1(\pi_1(\M);Ad \rho)$ unchanged.
The diffeomorphisms $f_t$ amount to choosing a different map
from $\M_0$ to $\M_t$. But $f_t$ is isotopic
to $f_0 = {\rm identity}$, so the lifts $\tilde f_t$
don't change the group cocycle at all.  It follows that equivalent families
of hyperbolic structures determine the same group cohomology class.

If, instead, we view the infinitesimal deformation as represented
by the $E$-valued $1$-form $\omega$, we note that the infinitesimal
effect of the isometries $k_t$ is to add a constant to
$s: \tilde \M \to \tilde E$.  Thus, $d s$, its
projection $\omega$, and the infinitesimal variation of metric are
all unchanged.  However, the infinitesimal
effect of the $\tilde f_t$ is to change the vector field on
$\tilde \M$ by the lift of a globally defined vector field
on $\M$.  This changes $\omega$ by the derivative of a
{\em globally defined} section of $E$.  Hence, it doesn't
affect the de Rham cohomology class in $H^1(\M; E)$.
The corresponding infinitesimal change of metric is altered
by the Lie derivative of a globally defined vector field
on $\M$.

\section{Infinitesimal harmonic deformations}\label{infharm}

In the previous section, we saw how a family of hyperbolic structures leads,
at the infinitesimal level, to both a group cohomology class and a
de Rham cohomology class.  Each of these objects has certain advantages and
disadvantages.  The group cohomology class is determined by its values on a
finite number of group generators and the equivalence relation, dividing out
by coboundaries which represent infinitesimal conjugation by a Lie 
group element,
is easy to understand.  Local changes in the geometry of the 
hyperbolic manifolds
are not encoded, but important global information like the 
infinitesimal change in
the lengths of geodesics is easily derivable from the group cohomology class.
However, the chosen generators of the fundamental group may not be related in
any simple manner to the hyperbolic structure, making it unclear how the
infinitesimal change in the holonomy representation affects the geometry of
the hyperbolic structure. Furthermore, it is usually hard to compute 
even the dimension
of $H^1(\pi_1(\M);Ad \rho)$ by purely algebraic means and much more 
difficult to find
explicit classes in this cohomology group.

The de Rham cohomology cocycle does contain information about the 
local changes in metric.
The value of the corresponding group cocycle applied to an element
$\gamma \in \pi_1(\M)$ can be
computed simply by integrating an $E$-valued $1$-form representing 
the de Rham class
around any loop in the homotopy class of $\gamma$ that element; this is the
definition of the de Rham isomorphism map.  However, it is generally 
quite difficult
to find such a $1$-form that is sufficiently explicit to carry out 
this computation.
Furthermore, the fact that any de Rham representative can be altered, 
within the
same cohomology class, by adding an exact $E$-valued $1$-form, (which 
can be induced by
any smooth vector field on $\M$), means that the behavior on small open sets is
virtually arbitrary, making it hard to extrapolate
to information on the global change in the hyperbolic metric.

In differential topology, one method for dealing with the large indeterminacy
within a {\em real-valued} cohomology class is to use Hodge theory. 
The existence
and uniqueness of a closed and co-closed (harmonic) $1$-form within a 
cohomology
class for a closed Riemannian manifold is now a standard fact.  Similar
results are known for complete manifolds and for manifolds with boundary,
where uniqueness requires certain asymptotic or boundary conditions 
on the forms.
By putting a natural metric on the fibers of the bundle $E$, the same theory 
extends to the de Rham cohomology groups, $H^1(\M; E)$, that arise in the 
deformation theory of hyperbolic structures.  The fact that these forms are 
harmonic implies that they satisfy certain nice elliptic linear partial 
differential equations.
In particular, for a harmonic representative $\omega \in H^1(\M; E)$, 
the infinitesimal
change  in metric $\eta$, which appears as the symmetric portion of 
the translational
part  of $\omega$, satisfies equations of this type.  As we will see 
in the next
section,  these are the key to the infinitesimal rigidity of 
hyperbolic structures.

For manifolds with hyperbolic metrics, the theory of harmonic maps provides
a non-linear generalization of this Hodge theory, at least for closed 
manifolds.
For non-compact manifolds or manifolds with boundary, the asymptotic 
or boundary
conditions needed for this theory are more complicated than those needed for
the Hodge theory.  However, at least the relationship described below 
between the
defining equations of the two theories continues to be valid in this 
general context.

It is known that, given
a map $f:\M \to \M^{\prime}$ between closed hyperbolic manifolds, there is
a unique harmonic map homotopic to $f$.  (In fact, this holds for
{\em negatively curved} manifolds. See \cite{ES}.)
Specifically, if $\M =\M^{\prime}$
and $f$ is homotopic to the identity, the identity map is this unique
harmonic map.  Associated to a $1$-parameter family $\M_t$ of hyperbolic
structures on $\M$ is a $1$-parameter family of developing maps from the
universal cover $\tilde \M$ of $\M$ to $\H^n$.  Using these maps to pull
back the metric on $\H^n$ defines a $1$-parameter family of metrics on
$\tilde \M$, and dividing out by the group of covering transformations
determines a family of hyperbolic metrics $g_t$ on $\M$.  However, a
hyperbolic structure only determines an {\em equivalence class} of developing
maps. Because of this equivalence relation, the metrics, $g_t$ are
only determined, for each fixed $t$, up to pull-back by a diffeomorphism of
$\M$.  For the smooth family of hyperbolic metrics $g_t$ on $\M$, we consider
the identity map as a map from $\M$, equipped with the metric $g_0$, to
$\M$, equipped with the metric $g_t$.  For $t=0$, the identity map is harmonic,
but in general it won't be harmonic.  Choosing the unique harmonic map
homotopic to the identity for each $t$ and using it to pull back the metric
$g_t$ defines a new family of metrics beginning with $g_0$.  (For small
values of $t$ the harmonic map will still be a diffeomorphism.)  By uniqueness
and the behavior of harmonic maps under composition with an isometry, the
new family of metrics depends only on the family of hyperbolic structures.
In this way, we can pick out a canonical family of metrics from the
family of equivalence classes of metrics.

If we differentiate this ``harmonic" family of metrics associated to
a family of hyperbolic structures at $t=0$,  we obtain a symmetric
$2$-tensor which describes the infinitesimal change of metric at each
point of $\M$.  Using the underlying hyperbolic metric, a symmetric $2$-tensor
on $\M$ can be viewed as a symmetric $T\M$-valued $1$-form.  This is 
precisely the
form $\eta$ described above which is the symmetric portion of the translational
part of the Hodge representative $\omega \in H^1(\M; E)$, corresponding to
this infinitesimal deformation of the hyperbolic structure.

Thus, the Hodge representative in the de Rham cohomology group corresponds
to an infinitesimal harmonic map.  The corresponding infinitesimal change
of metric has the property that it is $L^2$-orthogonal to the trivial
variations of the initial metric given by
the Lie derivative of compactly supported vector fields on $\M$.

We now specialize to the case of interest in this paper,
3-dimensional hyperbolic cone-manifolds.  We recall some of
the results and computations derived in \cite{HK1}.
Let $M_t$ be a smooth family of hyperbolic cone-manifold
structures on a 3-dimensional manifold $M$ with cone angles $\alpha_t$ along
a link $\Sigma$, where $0 \le \alpha_t  \le 2\pi$.
By the Hodge theorem proved in \cite{HK1},
the corresponding infinitesimal deformation at time $t=0$
has a unique Hodge representative whose translational part is a
$T\M$-valued 1-form $\eta$ on $\M = M-\Sigma$
satisfying
\begin{eqnarray} D^* \eta =0, \\
D^* D \eta = -\eta.
\label{harmeqs}\end{eqnarray}
Here $D: \Omega^1(\M;T\M) \to \Omega^2(\M;T\M)$ is the exterior
covariant derivative, defined, in terms of the Riemannian connection
from the hyperbolic metric on $\M$, by
$$D\eta(v,w)=\nabla_v\eta(w)-\nabla_w\eta(v)-\eta([v,w])$$
for all vectors fields $v,w$ on $\M$, and
$D^* : \Omega^2(\M;T\M) \to \Omega^1(\M;T\M)$ is its formal adjoint.
Further, $\eta$ and $*D\eta$ determine symmetric and traceless
linear maps $T_x \M \to T_x \M$ at each point $x \in \M$.

Inside an embedded tube $U=U_R$ of radius
$R$ around the singular locus $\Sigma$,
$\eta$ has a decomposition:
\begin{eqnarray}\eta = \eta_m + \eta_l + \eta_c \label{etadecomp}
\end{eqnarray}
where $\eta_m$, $\eta_l$ are ``standard''  forms
changing the holonomy of peripheral group elements,
and $\eta_c$ is a correction term with $\eta_c$,
$D\eta_c$ in $L^2$.

We think of $\eta_m $ as an ideal model for the infinitesimal
deformation in a tube around the
singular locus; it is completely determined by the
rate of change of cone angle. Its effect on
the complex length ${\cal L}$ of any peripheral element satisfies
$${d {\cal L} \over d \alpha} = {{\cal L} \over \alpha}.$$
In particular, the (real) length $\l$ of the core geodesic satisfies
\begin{eqnarray}\label{model}
{d \l \over d \alpha} = {\l \over \alpha}
\end{eqnarray}
for this model deformation.

This model is then ``corrected" by adding $\eta_l$ to get the actual change in
complex length of the core geodesic and
then by adding a further term $\eta_c$ that doesn't
change the holonomy of the peripheral elements at all, but is needed to
extend the deformation in the tube $U$ over the rest of the manifold $\M$.

One special feature of the 3-dimensional case is the
{\em complex structure} on the Lie algebra
$\g \cong sl_2 \C$ of infinitesimal isometries of $\H^3$.
The infinitesimal rotations fixing a point $p \in \H^3$
can be identified with $su(2) \cong so(3)$, and
the infinitesimal pure translations at $p$
correspond to $i \, su(2) \cong T_p \H^3$. Geometrically,
if $t \in T_p \H^3$ represents an infinitesimal translation,
then $i \, t $ represents an infinitesimal rotation with
axis in the direction of $t$. Thus,
on a hyperbolic 3-manifold $\M$ we can identify
the bundle $E$ of (germs of) infinitesimal isometries
with the {\em complexified} tangent bundle $T\M \otimes \C $.

In \cite{HK1} it was shown that the corresponding harmonic
$1$-form $\omega$ with values in the infinitesimal isometries of $\H^3$
can be written in this complex notation as:
\begin{eqnarray}
\omega = \eta + i *\! D\eta.
\label{omega}\end{eqnarray}
There is decomposition of $\omega$
in the neighborhood $U$ analogous to that (\ref{etadecomp}) of $\eta$ as
\begin{eqnarray}
\omega = \omega_m + \omega_l + \omega_c,\label{omegadecomp}
\end{eqnarray}
where only $\omega_m$ and $\omega_l$
change the peripheral holonomy and $\omega_c$ is in $L^2$.

The fact that the hyperbolic structure on $\M = M-\Sigma$ is incomplete
makes the existence and uniqueness of a Hodge representative substantially
more subtle than the standard theory for complete hyperbolic structures
(including structures on closed manifolds).  Certain conditions on the
behavior of the forms as they approach the singular locus are required.
This makes the theory sensitive to the value of the cone angle at the
singularity; in particular, this is where the condition that the cone
angle be at most $2\pi$ arises.  The fact that $\omega_c$ is in $L^2$
is a reflection of these asymptotic conditions.  In the final section of
this paper, we discuss a new version of this Hodge theory, involving
boundary conditions on the boundary of a tube around the singular locus,
that removes the cone angle condition, replacing it with a lower bound
on the radius of the tube.

\section{Effective Rigidity} \label{effrigid}
In this section we explain how the equations satisfied by
the harmonic representative of an infinitesimal deformation lead
to local rigidity results.  We then come to one of our primary
themes, that the arguments leading to local rigidity can be made
computationally effective.  By this we mean that even when there
{\it does} exist a non-trivial deformation of a hyperbolic structure,
the same equations can be used to bound both the geometric and analytic
effect of such a deformation.  This philosophy carries over into
many different contexts, but here we will continue to focus on
finite volume $3$-dimensional hyperbolic cone-manifolds.

The first step is to represent an infinitesimal deformation by a
Hodge (harmonic) representative $\omega$ in the cohomology group 
$H^1(\M;E)$, as
discussed in the previous section.  If $\M$ is any hyperbolic $3$-manifold,
the symmetric real part of this representative is a 1-form $\eta$ with values
in the tangent bundle of $\M$, satisfying the Weitzenb\"ock-type formula:
$$D^*D\eta ~+~ \eta ~=~ 0$$
where $D$ is the exterior covariant derivative on such forms and
$D^*$ is its adjoint.  First, suppose $\M$ is closed.  Taking the
$L^2$ inner product of this
formula with $\eta$ and integrating by parts gives the formula
$$||D\eta||_X^2 ~+~ ||\eta||_X^2 ~=~ 0.$$
(Here $||\eta||_X^2$ denotes the square of the
$L^2$ norm of $\eta$ on $X$.)
Thus $\eta ~=~ 0$ and the deformation is
trivial. This is the proof of local rigidity for
closed hyperbolic 3-manifolds, using the methods of
Calabi \cite{cal}, Weil \cite{weil} and Matsushima-Murakami \cite{MM}.

When $\M$ has boundary or is non-compact, there will be a
Weitzenb\"ock boundary term $\b$:
\begin{eqnarray}||D\eta||_X^2 ~+~ ||\eta||_X^2 ~=~ \b .\label{bdryeq}
\end{eqnarray}
If the boundary term is non-positive, the same conclusion
holds: the deformation is trivial.  When $\M = M - \Sigma$, where
$M$ is a hyperbolic cone-manifold with cone angles at most $2\pi$
along its singular set $\Sigma$, it was shown in \cite{HK1} that, for
a deformation which leaves the cone angle fixed, it is possible
to find a representative as above for which the boundary term
goes to zero on the boundary of tubes around the singular locus whose
radii go to zero.  Again, such an infinitesimal deformation must be trivial,
proving {\em local rigidity rel cone angles}.

On the other hand, Thurston has shown (\cite[Chap. 5]{thnotes}) that there
{\em are} non-trivial deformations of the (incomplete) hyperbolic structures
on $\M = M - \Sigma$.  By local rigidity rel cone angles such a 
deformation must
change the cone angles, implying that it is always possible to alter the cone
angles by a small amount.  Using the implicit function theorem it is further
possible to show that the variety of representations $\pi_1 (\M) \to 
PSL_2(\C)$ is
smooth near the holonomy representation of such a hyperbolic cone-manifold.
This leads to a local parametrization of hyperbolic cone-manifolds by 
cone angles.

\begin{theorem}[\cite{HK1}] For a $3$-dimensional hyperbolic cone-manifold with
singularities along
a link with cone angles $\le 2\pi$, there are no deformations
of the hyperbolic structure keeping the cone angles fixed.
Furthermore, the nearby hyperbolic cone-manifold structures are
parametrized by their cone angles.
\label{localparam}
\end{theorem}

The {\em argument} for local rigidity rel cone angles actually provides 
further information about the boundary term.  To explain this, we need
to give a more detailed description of some of the work in \cite{HK1}. 
This will provide not only a fuller explanation of the proof that there are
no deformations fixing the cone angles, but also additional information about 
the deformations that {\em do} occur.

Assume that $\eta$ represents a non-trivial infinitesimal deformation.
Recall that,
inside a tube around the singular locus, $\eta$ can be decomposed as
$\eta = \eta_m + \eta_l + \eta_c$, where only $\eta_m$ changes the cone angle.
Leaving the cone angle unchanged is equivalent to the vanishing
of $\eta_m$.  As we shall see below, the boundary term
for $\eta_m$ by itself is positive.  Roughly speaking, $\eta_m$ contributes
positive quantities to the boundary term, while everything else gives
negative contributions. (There are also cross-terms which are easily handled.)
The condition that the entire boundary term be positive not only implies 
that the $\eta_m$ term {\em must} be non-zero (which is equivalent to local
rigidity rel cone angles), but also puts strong restrictions on
the $\eta_l$ and $\eta_c$ terms.  This is the underlying philosophy for
the estimates in this section.

In order to implement this idea, we need to
derive a formula for the boundary term in (\ref{bdryeq}) as an integral
over the boundary of $\M$.
For details we refer to \cite{HK1}.

Let $U_r$ denote a tubular neighborhood of radius $r$
about the singular locus of $M$ and let $\M = M- U_r$; it will always
be assumed that $r$ is small enough so that $U_r$ will be embedded.
Let $T_r$ denote the boundary torus of $U_r$,
oriented by the normal ${\bd\over \bd r}$, (which is the {\em inward} normal
for $\M$).
For any $T\M$-valued $1$-forms $\alpha, \beta$ we define
\begin{eqnarray}\b_r(\alpha,\beta) = \int_{T_r} *D\alpha \wedge \beta.
\label{b_defn}\end{eqnarray}

Note that in this integral, $*D\alpha \wedge \beta$ denotes
the real valued 2-form obtained
using the wedge product of the form parts, and the geometrically
defined inner product on the vector-valued parts of the
$T\M$-valued $1$-forms $*\!D\alpha$ and $\beta$.

As above, we express the Hodge $E$-valued 1-form as $\omega = \eta + 
i\, *\!D\eta$
where $D^{*}D \eta + \eta = 0$.  Fix a radius $R$, remove the tubular
neighborhood $U_R$, and denote $M- U_R$ by $\M$.
Then one computes that the Weitzenb\"ock boundary term $\b$ in 
(\ref{bdryeq}) equals
$\b_R(\eta,\eta)$ (see Proposition 1.3 and p. 36 of \cite{HK1}).
This implies:

\begin{lemma}
\begin{eqnarray}  \b_R(\eta,\eta) =  ||\eta||_{\M}^2 + ||D\eta||_{\M}^2
= ||\omega||_{\M}^2 .
\label{etabdry}\end{eqnarray}
\end{lemma}

In particular, we see that $b_R(\eta,\eta)$ is {\em non-negative}.
Writing $\eta = \eta_0 + \eta_c$ where $\eta_0 = \eta_m + \eta_l$,
we analyze the contribution from each part.
First, using the Fourier decomposition
for $\eta_c$ obtained in \cite{HK1}, it turns out that
the cross-terms vanish so that
the boundary term is simply the sum of two boundary terms:
\begin{eqnarray} \label{bdrydecomp} \b_R(\eta,\eta) = \b_R(\eta_0,\eta_0)
+ \b_R(\eta_c,\eta_c).
\end{eqnarray}

Next, we see that the contribution, $\b_{R}(\eta_c,\eta_c)$,
from the part of the ``correction term" that doesn't affect
the holonomy of the peripheral elements, is {\em non-positive}.  In fact,

\begin{proposition}\label{etacbdryprop}
\begin{eqnarray} \b_{R}(\eta_c,\eta_c) =
- ( ||\eta_c||_{U_R}^2 + ||D\eta_c||_{U_R}^2) =
- ||\omega_c||_{U_R}^2.
\label{etacbdry}
\end{eqnarray}
\end{proposition}

We have assumed that
$\omega_c$ is harmonic in a neighborhood of $U_R$ so the
same argument applied above to $\eta$ can be applied to $\eta_c$ on
this neighborhood.  Consider a region $N$ between tori at distances
$r, R$ from $\Sigma$ with $r < R$.  As before, integration by parts over
this region implies that the difference $\b_r(\eta_c,\eta_c) - 
\b_R(\eta_c,\eta_c)$
equals $||\omega_c||_N^2$.
Then the main step is to show that
$\lim_{r\to 0} \b_r(\eta_c,\eta_c) = 0$.
This follows from the
proof of rigidity rel cone angles (in section 3 of \cite{HK1}),
since $\eta_c$ represents an infinitesimal deformation which
doesn't change the cone angle.

Combining (\ref{bdrydecomp}) with (\ref{etabdry}) and
(\ref{etacbdry}), we obtain:
\begin{eqnarray}
\b_R(\eta_0,\eta_0) = ||\omega||_{M-U_R}^2 + ||\omega_c||_{U_R}^2.
\label{eta0bdry}\end{eqnarray}

In particular, this shows that
\begin{eqnarray} 0 \leq \b_R(\eta_0,\eta_0), \label{eta0pos}\end{eqnarray}

\smallskip

and that
\begin{eqnarray}||\omega||_{M-U_R}^2 \leq \b_R(\eta_0,\eta_0).\label{normbound}
\end{eqnarray}

\begin{remark} We emphasize that the only place in the derivation of
(\ref {eta0pos}) and (\ref{normbound}) that we have used the analysis
near the singular locus from \cite{HK1} is in the proof of 
Proposition \ref{etacbdryprop}.
Furthermore, all that is required from this Proposition is the fact that
$\b_{R}(\eta_c,\eta_c)$ is {\em non-positive}.  This, together with
(\ref{etabdry}) and (\ref{bdrydecomp}), implies both (\ref{eta0pos}) 
and (\ref{normbound}).
In the final section of this paper, we describe another method for finding
a Hodge representative for which $\b_{R}(\eta_c,\eta_c)$ is non-positive.
This method requires a tube radius of at least a universal size, but 
no bound on
the cone angle.  Once this is established, all the results described here
carry over immediately to the case where the tube radius condition is 
satisfied.
\end{remark}

We will focus here on applications of the inequality (\ref{eta0pos})
which is the primary use of this analysis in \cite{HK2}.  The work of
Brock and Bromberg discussed in these proceedings (\cite{BB2}) also
requires the second inequality (\ref{normbound}).   This is discussed
further in the final section.

As we show below, (\ref{eta0pos}) implies that, in the decomposition
$\eta_0 = \eta_m + \eta_l$, the $\eta_m$ term must be non-trivial
for a non-trivial deformation.  This is equivalent to local rigidity 
rel cone angles.
The positivity result (\ref{eta0pos}) can also
be used to find {\em upper bounds} on $\b_R(\eta_0,\eta_0)$.
On the face of it, this may seem somewhat surprising, but, as we explain
below, the algebraic structure of the quadratic form $\b_R(\eta_0,\eta_0)$
makes it quite straightforward to derive such bounds.

The possible harmonic forms $\eta_0= \eta_m + \eta_l$ give a
3-dimensional real vector space $W$ representing models for deformations
of hyperbolic cone-manifold structures in a neighborhood of the boundary torus.
Here $\eta_m$ lies in a 1-dimensional subspace $W_m$ containing
deformations changing the cone angle, while $\eta_l$
lies in the 2-dimensional subspace $W_l$ consisting of deformations
leaving the cone angle unchanged. In \cite{HK1}, we describe
explicit $TX$-valued 1-forms giving bases for these subspaces.

Now consider the quadratic form $Q(\eta_0) = b_R(\eta_0,\eta_0)$ on 
the vector space $W$.
One easily computes that in all cases, $Q$ is {\em positive definite} on $W_m$
and {\em negative definite} on $W_l$, so $Q$ has signature $+--$. 
This gives the
situation shown in the following figure.

\medskip
\centerline{\epsfbox{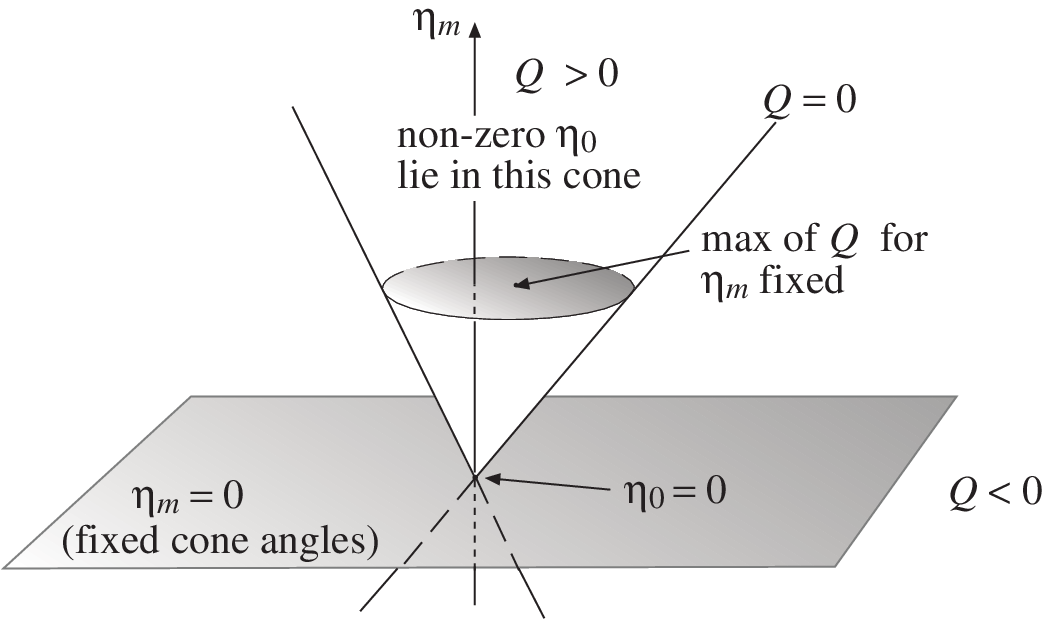}}
\medskip

The positivity condition (\ref{eta0pos}) says that
$\eta_0$ lies in the cone where $Q\ge 0$ for any deformation which 
extends over the manifold
$\M$.
Further, $\eta_m$ must be non-zero if the deformation is non-trivial;
so the cone angle must be changed. Thus (\ref{eta0pos}) implies
local rigidity rel cone angles.

As noted above, the local parametrization by cone angles
(Theorem \ref{localparam}) follows from this, and a
smooth family of cone-manifold structures $M_t$ is completely
determined by a choice of parametrization of the cone angles $\alpha_t$.
We are free to choose this parametrization as we wish.
Then the term $\eta_m$ is completely determined by
the derivative of the cone angle.

Once $\eta_m$ is fixed, the inequality $Q(\eta_0) \ge 0$ restricts
$\eta_0 =\eta_m + \eta_l$ to lie in an ellipse (as illustrated above).
Since $Q$ has a positive maximum
on this compact set, this
gives an explicit upper bound for $Q(\eta_0)$ for any deformation.

\medskip
It turns out to be useful to parametrize the cone-manifolds by the
{\em square} of the cone angle $\alpha$; i.e., we will let $t = \alpha^2$.
With this choice of parametrization we obtain:
 
\begin{eqnarray}
\b_R(\eta_m,\eta_m) ~=~ {\area(T_R) \over 16 m^4}~ ({\rm tanh}\,R + 
{\rm tanh}^3\, R) ,
\label{bmm}\end{eqnarray}
where $m =  \alpha \sinh(R)$ is the length of the meridian on the
tube boundary $T_R$.

\smallskip

Essentially, any contribution to $\b_R(\eta_0,\eta_0)$ from
$\eta_l$ will be negative; the cross-terms only complicate matters slightly.
Computation leads easily to the following upper bound:

\begin{eqnarray}
\b_R(\eta_0,\eta_0) \le {\area(T_R) \over 8 m^4}.
\label{b00}\end{eqnarray}

\smallskip

Combining (\ref{eta0bdry}) and (\ref{b00}), we see that
the boundary formula (\ref{bdryeq}) leads to the estimate:

\begin{theorem}
\begin{eqnarray}
||\omega||_{M-U_R}^2 + ||\omega_c||_{U_R}^2
\le {\area(T_R) \over 8 m^4},
\label{L2estimate}\end{eqnarray}
where $m =  \alpha \sinh(R)$ is the length of the meridian on the
tube boundary $T_R$.
\end{theorem}

\begin{remark}
A crucial property of the inequalities (\ref{b00}) and (\ref{L2estimate})
is their dependence only on the geometry of the boundary torus, not
on the rest of the hyperbolic manifold.  This is the reason that
{\em a priori} bounds, independent of the underlying manifold can be
derived by these methods.  Furthermore, it is possible to find geometric
conditions, like a very short core geodesic, that force the upper bounds
to be very small.

In particular, (\ref{L2estimate}) provides an upper bound on
the $L^2$ norm of $\omega$ on the complement of the tubular neighborhood
of the singular locus.  Such a bound can be used to bound the infinitesimal
change in geometric quantities, like lengths of geodesics, away from
the singular locus.  Arbitrarily short core geodesics typically lead to the
conclusion that these infinitesimal changes are arbitrarily small, a useful
fact when studying approximation by manifolds with short geodesics.
This will be discussed further in the final section.
\end{remark}

The principle that the contribution of $\eta_l$ to $\b_R(\eta_0,\eta_0)$
is essentially negative (ignoring the cross-terms) while
$\b_R(\eta_0,\eta_0)$ is positive also provides a bound on the size of
$\eta_l$.
Although $\eta_m$ does change the length of the singular locus,
its effect is fixed since $\eta_m$ is fixed by the parametrization.
The only other term affecting the length of the singular locus is $\eta_l$.
Explicit computation ultimately leads to the following estimate:

\begin{theorem}[\cite{HK2}] \label{keyboundtheorem}
Consider any smooth family of hyperbolic cone structures on $M$, all
of whose cone
angles are at most $2 \pi$.  For any component of the singular set, let
$\l$ denote its length and $\alpha$ its cone angle.  Suppose there is
an embedded tube of radius $R$ around that component.  Then
$${d \l \over d \alpha} = {\l \over \alpha} ( 1 + E ),$$
where
\begin{eqnarray} {-1 \over \sinh^2(R)}
	\biggl({ 2\sinh^2(R) + 1 \over 2 \sinh^2(R) + 3}\biggr)
   \le E \le {1\over \sinh^2(R) }.
\label{keybound}\end{eqnarray}
\end{theorem}

Note that the `error' term $E$ represents the deviation from the
standard model; compare (\ref{model}).

\smallskip

This is the key estimate in \cite{HK2}.  The
next section discusses some of the applications of this inequality
to the theory of hyperbolic Dehn surgery.

\section{A quantitative hyperbolic Dehn surgery theorem} \label{dehn}

We begin with a non-compact, finite volume hyperbolic $3$-manifold
$\M$, which, for simplicity, we assume has a single cusp.  In the general
case the cusps are handled independently.  The manifold $\M$ is the
interior of a compact manifold which has a single torus boundary, which
is necessarily incompressible.  By attaching a solid torus by a diffeomorphism
along this boundary torus, one obtains a closed manifold, determined
up to diffeomorphism by the isotopy class of the non-trivial simple
closed curve $\gamma$ on the torus which bounds a disk in the solid torus.
The resulting manifold is denoted by $\M(\gamma)$ and $\gamma$ is referred to
either as the {\em meridian} or the {\em surgery curve}.  The process is called
{\em Dehn filling} (or Dehn surgery).

The set of Dehn fillings of $\M$ is thus parametrized by the set of
simple closed curves on the boundary torus; after choosing a basis
for $H_1(T^2,\Z)$, these are parametrized by pairs $(p,q)$ of relatively
prime integers.  Thurston (\cite{thnotes}) proved that, for all but
a finite number of choices of $\gamma$, $\M(\gamma)$ has a complete,
smooth hyperbolic structure.  However, the proof is computationally
ineffective.  It gives no indication of how many non-hyperbolic fillings
there are or which curves $\gamma$ they might correspond to.  In
particular, it left open the possibility that there is no upper bound to the
number of non-hyperbolic fillings as one varies over all possible $\M$.

One approach to putting a hyperbolic metric on $\M(\gamma)$ is through
families of hyperbolic cone-manifolds.
The complete metric on the open manifold $\M$ is deformed
through incomplete metrics whose metric completions are hyperbolic
cone-manifold structures on $\M(\gamma)$, with the singular set equal to
the core of the added solid torus.  The complete structure can be considered
as a cone-manifold with angle $0$.  The cone angle is increased
monotonically, and, if the angle of $2 \pi$ is reached, it defines
a smooth hyperbolic metric on $\M(\gamma)$.

Thurston's proof of his finiteness theorem actually shows that,
for {\em any} non-trivial simple closed curve $\gamma$ on the boundary torus,
it is possible to find such a cone-manifold structure on $\M(\gamma)$
for sufficiently small (possibly depending on $\M$) values of the cone angle.
This also follows from Theorem \ref{localparam}, which further 
implies that, for
any angle at most $2 \pi$,  it is always possible to
change the cone angle a small amount, either to increase it or to decrease
it.  Locally, this can be done in a unique way since
the cone angles locally parametrize the set of cone-manifold structures
on $\M(\gamma)$.
Although there are always small variations of the cone-manifold structure,
the structures may degenerate in various ways as a family of
angles reaches a limit.  In order to find a smooth hyperbolic metric
on $\M(\gamma)$ it is necessary and sufficient to show that no 
degeneration occurs
before the angle $2 \pi$ is attained.

In \cite{HK2} the concept of convergence of metric spaces in the 
Gromov-Hausdorff
topology (see, e.g., \cite{GrFrench}, \cite{GrEnglish}, or \cite{CHK1})
is utilized to rule out degeneration of the metric as long as there
is a lower bound on the tube radius and an upper bound on the volume.
The main issue is to show that the injectivity radius in the complement
of the tube around the singular locus stays bounded below.  One shows that,
if the injectivity radius goes to zero, a new cusp develops.  Analysis of
Dehn filling on this new cusp leads to a contradiction of Thurston's
finiteness theorem.  The argument is similar to a central step in the proof of
the orbifold theorem (see, e.g., \cite{CHK1}, \cite{BP}, or \cite{BLP}).
One result proved in \cite{HK2} is:

\begin{theorem} \label{thm2}  Let $\N_t$, $t\in [0,t_\infty)$, be
a smooth path of closed hyperbolic cone-manifold structures
on $(\N, \Sigma)$ with cone angle $\alpha_t$ along the singular locus $\Sigma$.
Suppose $\alpha_t \to \alpha \geq 0$ as $t \to t_\infty$, that the volumes of
the $\N_t$ are bounded above, and that there is a positive constant
$R_0$ such that there is an embedded tube of radius at least $R_0$ around
$\Sigma$ for all $t$.  Then the path extends continuously to $t=t_\infty$
so that as $t \to t_\infty$, $\N_t$ converges in the bilipschitz topology
to a cone-manifold structure $\N_\infty$ on $\N$ with cone angles
$\alpha$ along $\Sigma$.
\end{theorem}

As a result of this theorem, we can focus on controlling the tube radius.
A general principle for smooth hyperbolic manifolds is that short geodesics
have large embedded tubes around them.  This follows from the Margulis lemma
or, alternatively in dimension $3$, from J{\o}rgensen's inequality. 
However, both of
these results require that the holonomy representation be discrete which
is almost never true for cone-manifolds.  In fact, the statement that there
is a universal lower bound to the tube radius around a short geodesic is easily
seen to be false for cone-manifolds.  To see this, consider a sequence of
hyperbolic cone-manifolds whose diameters go to zero.  Then, both the
length of the singular locus and the tube radius go to zero.  For 
example, $S^3$ with
singular locus the figure-$8$ knot and cone angles approaching ${2 
\pi \over 3}$
from below behaves in this manner.

However, a more subtle statement is true.  If, at the beginning of
a family of hyperbolic cone-manifold structures, the tube radius is
sufficiently large and the length of the core curve is sufficiently small,
then as long as the core curve remains sufficiently small, the tube
radius will be bounded from below.  The precise statement is given below.
It should be noted that it is actually the {\em product} of the cone angle
and the core length that must be bounded from above.

\begin{theorem}[\cite{HK2}] \label{alnondeg} Let $\N_s$ be a smooth
family of finite volume $3$-dimensional hyperbolic cone-manifolds, with
cone angles $\alpha_s, 0\le s < 1$, where $\lim_{s\to 1}\alpha_s = \alpha_1$.
Suppose the tube radius $R$ satisfies $R \ge 0.531$ for $s=0$
and $\alpha_s \l_s \le 1.019$ holds for all $s$,
where $\l_s$ denotes the
length of the singular geodesic.  Then the tube radius satisfies
$R \ge 0.531$ for all $s$.
\end{theorem}

The proof of this theorem involves estimates on ``tube packing"
in cone-manifolds.  We look in a certain cover of the complement
of the singular locus where the tube around the singular geodesic
lifts to many copies of a tube of the same radius around infinite
geodesics which are lifts of the core geodesic.  One of the
lifts is chosen and then packing arguments are developed to show that
other lifts project onto the chosen one in a manner that fills up
at least a certain amount of area in the tube in the original cone-manifold.
In this way, we derive a lower bound for the area of the boundary torus which
depends only on the product $\alpha \l$, assuming $R \ge 0.531$.
This, in turn, further bounds the tube radius {\em strictly} away
from $0.531$ as long as $\alpha \l$ is sufficiently small.  Thus,
if this product stays small the tube radius stays away from the value
$0.531$ so that the estimates continue to hold.  The estimates derived
actually prove that, as long as these bounds hold, then $R \to \infty$
as $\alpha \l \to 0$.  Furthermore, the rate at which the tube radius
goes to infinity is bounded below.

This result shows that the tube radius can be bounded below by
controlling the behavior of the holonomy of the peripheral elements.
The cone angle is determined by the parametrization so it suffices to
understand the longitudinal holonomy.  In the previous section, we
derived estimates (\ref{keybound}) on the derivative of the core length
with respect to cone angle, where the bounds depend on the tube radius.
On the other hand, from Theorem \ref{alnondeg}, the change in the tube radius
can be controlled by the product of the cone angle and the core length.
Putting these results together, we arrive at differential
inequalities that provide strong control on the change in the
geometry of the maximal tube around the singular geodesic, including
the tube radius.

A horospherical torus which is a cross-section of the cusp for the complete
structure on $\M$ has an intrinsic flat structure (i.e., zero 
curvature metric).
Any two such cross-sections differ only by scaling.  Given a choice $\gamma$
of surgery curve, an important quantity associated to this flat structure
is the {\em normalized flat length} of $\gamma$, which, by definition
equals the geodesic length of $\gamma$ in the flat metric, scaled to have area
$1$.  Clearly this is independent of the choice of cross-section.  Its
significance comes from the fact that its square (usually called the
{\em extremal length} of $\gamma$) is the limiting value of the ratio
${\alpha \over \l}$ of the cone angle to the core length
as $\alpha \to 0$ for hyperbolic cone-manifold structures on $\M(\gamma)$.
In particular, near the complete structure, even though
$\l$ and $\alpha$ individually approach zero, their {\it ratio}
approaches a finite, non-zero value.

The estimate (\ref{keybound}) implies that, as long as the tube radius
isn't too small, ${d\l \over d \alpha}$ is approximately equal to the
ratio ${\l \over \alpha}$.  In the case of equality (when the error term $E$
in (\ref{keybound}) equals zero), the ratio ${\l \over \alpha}$ will 
remain constant.
A small error term implies that the ratio doesn't change too much.
If the initial value of the {\it reciprocal}, ${\alpha \over \l}$ is large,
then ${\l \over \alpha}$ will be small
and stay small as long as the tube radius doesn't get too small.
But this implies that the product of the cone angle
and the core length will remain small.  In turn, the packing argument then
provides a lower bound to the tube radius.

Formally, this can be expressed as a differential inequality, bounding the
change of core length with respect to the change in cone angle as a function
of the core length and the cone angle.  Solving this inequality, with initial
conditions coming from the normalized length, gives a proof
of the following theorem:

\begin{theorem}[\cite{HK2}] \label{thm1}  Let $\M$ be a complete, finite volume
orientable hyperbolic $3$-manifold with one cusp and
let $T$ be a horospherical torus embedded as a cusp cross section. Fix
$\gamma$, a simple closed curve on $T$.
Let $\M_\alpha(\gamma)$ denote the
cone-manifold structure on $\M(\gamma)$ with cone angle $\alpha$
along the core, $\Sigma$, of the added solid torus, obtained by increasing
the angle from the complete structure.
If the normalized length of $\gamma$ on $T$ is at least $7.515$,
then there is a positive lower bound to the tube radius around $\Sigma$
in $\M_\alpha(\gamma)$ for all $\alpha$ satisfying $2\pi \ge \alpha \ge 0$.
\end{theorem}

Given $\M$ and $T$ as in Theorem \ref{thm1}, choose {\em any}
non-trivial simple closed curve $\gamma$ on $T$.  There is a
maximal sub-interval $J \subseteq [0, 2\pi]$ containing $0$ such
that there is a smooth family $\N_\alpha$, where $\alpha \in J$,
of hyperbolic cone-manifold structures on $\M(\gamma)$ with
cone angle $\alpha$.  Thurston's Dehn surgery theorem (\cite{thnotes})
implies that $J$ is non-empty and Theorem \ref{localparam}
implies that it is open.  Theorem \ref{thm2} implies that, with a lower bound
on the tube radii and an upper bound on the volume, the path of
$\N_\alpha$'s can be extended continuously to the endpoint of $J$.
Again, Theorem \ref{localparam}
implies that this extension can be made to be smooth.  Hence,
under these conditions $J$ will be closed.  By Schl\"afli's formula
(see (\ref{schlafli}) below), the volumes decrease as the cone angles increase
so they will clearly be bounded above.  Theorem \ref{thm1} provides
initial conditions on $\gamma$ which guarantee that there will be
a lower bound on the tube radii for all $\alpha \in J$.  Thus,
assuming Theorems
\ref{thm1} and \ref{thm2}, we have proved:

\begin{theorem}  \label{thm3}
Let $\M$ be a complete, finite volume orientable hyperbolic $3$-manifold with
one cusp, and let $T$ be a horospherical torus which is embedded as
a cross-section to the cusp of $\M$. Let
$\gamma$ be a simple closed curve on $T$ whose
normalized length is at least $7.515$.
Then the closed manifold $\M(\gamma)$ obtained by Dehn filling along
$\gamma$ is hyperbolic.
\end{theorem}

This result also gives a universal bound on the
number of non-hyperbolic Dehn fillings on a cusped hyperbolic $3$-manifold
$\M$, independent of $\M$.

\begin{corollary} \label{thm4}
Let $\M$ be a complete, orientable hyperbolic $3$-manifold with
one cusp. Then at most $60$ Dehn fillings on $\M$ yield manifolds
which admit no
complete hyperbolic metric.
\end{corollary}

When there are multiple cusps the results are only slightly weaker.
Theorem \ref{thm2} holds without change. If there are $k$ cusps,
the cone angles $\alpha_t$ and $\alpha$ are simply interpreted as
$k$-tuples of angles.  Having tube radius at least $R$ is interpreted
as meaning that there are disjoint, embedded tubes of radius $R$
around all of components of the singular locus.  The conclusion of
Theorem \ref{thm1} and hence of Theorem \ref{thm3} holds when there
are multiple cusps as long as the normalized lengths of all of the
meridian curves are at least $7.515 \sqrt{2}\approx 10.627$.
At most $114$ curves from each cusp need to be excluded.  In fact,
this can be refined to say that at most $60$ curves need to be
excluded from one cusp and at most $114$ excluded from the remaining
cusps.  The rest of the Dehn filled manifolds are hyperbolic.

By Mostow rigidity the volume of a closed or cusped hyperbolic $3$-manifold is
a topological invariant.  The set of volumes of hyperbolic $3$-manifolds
is well-ordered (\cite {thnotes}); the hyperbolic volume gives
an important measure
of the complexity of the manifold.  It is therefore of interest to 
find the smallest
volume hyperbolic manifold.  This is conjectured to be the Weeks 
manifold which can be
described as a surgery on the Whitehead link and has volume $\approx 0.9427$.
On the other hand, the smallest volume of a hyperbolic $3$-manifold 
with a single
cusp is known (\cite{CM}) to equal $\approx 2.0299$, the volume of the figure
eight knot complement.

Thus, one can attempt to estimate the volume of a closed hyperbolic
$3$-manifold by comparing it to the complete structure on the non-compact
manifold obtained by removing a simple closed geodesic from the 
closed manifold.
It is conjectured that these two hyperbolic structures are always connected
by a smooth family of cone-manifolds, with singular locus equal to the simple
closed geodesic, with cone angle $\alpha$ varying from $2 \pi$ in the smooth
structure on the closed manifold to $0$ for the complete structure on the
complement of the geodesic.   In such a family, Schl\"afli's formula 
implies that
the derivative of the volume $V$ satisfies the equation (see, e.g. 
\cite{Ho} or \cite[Theorem
3.20]{CHK1}):
\begin{eqnarray} {d V \over d\alpha} = -{1\over 2} \l,
\label{schlafli}\end{eqnarray}
\noindent where $\l$ denotes the length of the singular locus.  Thus, 
controlling of the
length of the singular locus throughout the family of cone-manifolds would
control the change in the volume.  In \cite{HK2} it is shown that
the derivative of the length of the singular locus with respect to
the cone angle is positive as long as the tube radius is at least
$\arctanh(1/\sqrt{3})$ and much sharper statements are proved,
using the packing arguments, when the length of the singular locus is
small.

It is also shown that, if the original simple closed geodesic
is sufficiently short, then such a family of cone-manifolds connecting
the smooth structure to the complete structure on the complement of the
geodesic will always exist.  To see this, note that, for sufficiently
short geodesics in a smooth structure, there will always be a tube
of radius greater than $\arctanh(1/\sqrt{3})$.  Thus the core length
will decrease as the cone angle decreases and, in particular, the product
of the cone angle and the core length will decrease.  By Theorems 
\ref{alnondeg}
and \ref {thm2} there can be no degeneration and the complete structure will
be reached.  (It is not hard to show that the volume is bounded above 
during the
deformation.)

Combining these two ideas, we can bound the volume of a closed
hyperbolic $3$-manifold with a sufficiently short geodesic in terms of
the associated cusped $3$-manifold.  An example of an explicit
estimate derived in this manner in \cite{HK2} is:

\begin{theorem}
Let $M$ be a {\em closed} hyperbolic manifold whose shortest closed
geodesic $\tau$ has length at most $0.162$.
Then the hyperbolic structure on $M$ can
be deformed to a complete hyperbolic structure on $M-\tau$ by decreasing
the cone angle along $\tau$ from $2\pi$ to $0$.
Furthermore, the volumes of these manifolds satisfy the inequality
${\rm Volume}(M) \ge {\rm Volume}(M - \tau) - 0.329.$
In particular,
${\rm Volume}(M) \ge 1.701$ so that it has larger volume than the
closed hyperbolic manifold with the smallest known volume.
\label{connect}
\end{theorem}

\section{Kleinian groups and boundary value theory}\label{boundaryvalues}

In this section we give a brief description of generalizations and applications
of the deformation theory described in the previous sections.

The first way in which one could attempt to generalize the harmonic
deformation theory is to allow finitely many infinite volume, but 
geometrically finite,
ends in our hyperbolic cone-manifolds $\M$.  Without giving a formal 
definition,
this structure provides ends that, asymptotically, are like the ends of
smooth geometrically finite hyperbolic $3$-manifolds.  We further 
assume that there are
no rank-$1$ cusps (so the ends are like those  Kleinian groups with 
{\em compact} convex
cores).  The term ``geometrically  finite" will include this extra 
assumption throughout
this section.   Each infinite volume end of such a geometrically 
finite cone-manifold
determines a conformal structure (in fact a complex projective structure)
on a surface at infinity, which can be used to (topologically) compactify
the cone-manifold.  We refer to these as {\em boundary conformal structures}.

When there is no singular set, then the quasi-conformal deformation 
theory developed
by Ahlfors, Bers and others (see, e.g., \cite{Ahl}, \cite{Bers})
implies that such structures are parametrized
by these conformal structures at infinity.  In particular, they satisfy
local rigidity relative to the boundary conformal structures; it is 
not possible
locally to vary the hyperbolic structure without varying at least one of the
conformal structures at infinity.

In \cite{Br1} Bromberg extends the harmonic deformation theory outlined
in the previous sections to such geometrically finite cone-manifolds, 
assuming that
the cone angles are at most $2 \pi$.  In particular, it is proved there that
there are no infinitesimal deformations of such structures that fix 
both the cone
angles and the conformal structures at infinity, i.e., rigidity rel
cone angles and boundary conformal structures.  This generalizes
the local theory both for smooth geometrically finite Kleinian groups and
for finite volume cone-manifolds.  However, it should be pointed out that,
since the holonomy groups for these geometrically finite cone-manifold
are not usually discrete, the global quasi-conformal theory on the
sphere at infinity which is the basis for the smooth (Ahlfors-Bers) theory
can't be used at all.  So, new techniques are required.  Similarly, since the
deformation is only assumed to be ``conformal at infinity", not 
trivial at infinity,
one cannot assume that there is a representative of this 
infinitesimal deformation for
which the infinitesimal change of metric is asymptotically trivial. 
As a result,
the asymptotic behavior of the $E$-valued forms representing the deformation
must be analyzed at these infinite volume ends as well as near the 
singular locus.
This makes the proof of the necessary Hodge theorem much more difficult.
After the Hodge theorem is proved, the final step is to show that there is a
harmonic (Hodge) representative for which the contribution to the Weitzenb\"ock
boundary term goes to zero near infinity.  This requires some subtle analysis.

Once this is established, the analytic results from the finite
volume cone-manifold theory go through without change.  In particular,
the inequalities (\ref {eta0pos}) and (\ref {normbound}) still apply.
Once packing arguments, which imply non-degeneration results, are proved in
this context, then the Dehn filling results, such as Theorem \ref{thm3},
can be generalized (with different numerical values).
Similarly, there will be an analog, for smooth, geometrically finite
hyperbolic manifolds, of Theorem \ref{connect}.  Recall that this theorem
says that a closed, smooth hyperbolic manifold with a sufficiently short
simple closed geodesic can be connected by a family of cone-manifolds to
the complete structure on the manifold obtained by removing the geodesic.

Such packing and non-degeneration results are proved by Bromberg in
\cite{Br3}.  However, his goal and the goal of subsequent papers
by (various subsets of) Brock, Bromberg, Evans and Souto
(see for example, \cite{Br2}, \cite{BB1}, \cite{BBES}) is an analytically
sharper version of this type of result.  A goal in each of these papers is
to approximate structures with very short (possibly singular) geodesics
by ones with smaller cone angles, including ones with cusps.
Not only is a path of cone-manifolds connecting the structures needed
but this must be constructed so as to bound the distortion
of the structure along the way.

As was remarked in Section \ref{effrigid}, the results in \cite{HK2}
which we have discussed in the previous sections require only the
inequality (\ref{eta0pos}).  However, the analysis in that section
also led to the inequality (\ref {normbound}) which bounds the
$L^2$ norm of the harmonic representative for the infinitesimal deformation.
Further calculations there also provided an upper bound 
(\ref{L2estimate}) for this
$L^2$ norm in terms of the geometry of the boundary torus.  It is easy to check
that, for any given non-zero cone angle, this upper bound will become 
arbitrarily
small as the length of the singular locus goes to zero, assuming there is
a lower bound to the tube radius.

This allows one to bound the $L^2$ norm of the harmonic representative
$\omega$ of an infinitesimal deformation at such a structure, where we
are assuming that the deformation is infinitesimally conformal
at infinity on the infinite volume ends.  Such deformations are 
locally parametrized
by cone angles.  To be consistent with the previous sections, we use the
parametrization $t = \alpha^2$, where $\alpha$ is the cone angle.
Then the following theorem holds:

\begin{theorem}[\cite{Br3}, \cite{HK3}] Let $\M$ be a geometrically finite
hyperbolic cone-manifold with no rank-$1$ cusps and
with an embedded tube of radius at least $\arctanh({1/\sqrt{3}})$ around the
singular locus $\tau$.   Suppose that $\omega$ is the harmonic 
representative of an
infinitesimal deformation of $\M$ as above.  Then,
for any fixed cone angle $\alpha > 0$ and any $\epsilon > 0$, there
is a length $\delta > 0$, depending only on $\alpha$ and $\epsilon$,
so that, if the length of $\tau$ is less than $\delta$, then
$||\omega||^2_{\M - U} < \epsilon$, for some embedded tube $U$
around $\tau$.
\label{shortgeo}
\end{theorem}

Recall that the real part of $\omega$ corresponds to the infinitesimal
change in metric induced by the infinitesimal deformation.  So, in particular,
the above theorem gives an $L^2$ bound on the size of the 
infinitesimal change of
metric.  However, it is still necessary to bound the change in the hyperbolic
structure in a more usable way.  In \cite{Br2} and \cite{Br3} this 
was turned into
a bound in the change of the projective structures at infinity.  This 
is sufficient
for the applications in those papers.  For the applications in \cite{BB1} and
\cite{BBES}, an upper bound is needed on the bilipschitz constants of 
maps between
structures along the path of cone-manifolds.  For such a bound, it is 
necessary to
turn the $L^2$ bounds on the infinitesimal change of metric into 
pointwise bounds.
It is then also necessary to extend the bilipschitz maps into the 
tubes in a way
that still has small bilipschitz constant.  This is carried out in \cite{BB1}.

The work in \cite{Br3} and \cite{BB1} globalizes Theorem \ref{shortgeo}.
It implies that, under the same assumptions, it is possible to find a
path of cone-manifolds from the geometrically finite cone-manifold $\M$ to the
complete structure on $\M$ with the geodesic removed, and that this can be done
so that all the geometric structures along the path can be made arbitrarily
close to each other.  In \cite{Br3} the distance between the structures is
measured in terms of projective structures at infinity, whereas in \cite{BB1}
it is measured by the bilipschitz constant of maps.  Since the
core geodesic is removed in the cusped structure, the authors call these
results ``Drilling Theorems'' (see \cite{BB2}).
These theorems provide very strong quantitative
statements of the qualitative idea that hyperbolic structures with short
geodesics are ``close" to ones with cusps.

The careful reader will have noticed that Theorem \ref{shortgeo} has
no conditions on the size of the cone angles whereas the theory
in \cite{HK1} requires that the cone angles be at most $2 \pi$.
As stated, this theorem and hence the full Drilling Theorem depends on a
harmonic deformation theory which has no conditions on the cone angle,
but only requires the above lower bound on the tube radius. Such a
theory is developed in \cite {HK3}.

Some uses of the Drilling Theorem (e.g., \cite{BBES}) only involve going
from a smooth structure (cone angle $2 \pi$) to a cusp, so only the analysis
in \cite{HK1} is needed.  However, the proofs of the Density Conjecture
in \cite{Br2}, \cite{Br3}, and \cite{BB1}, as described in these
proceedings (\cite{BB2}), require a path of cone-manifolds beginning with
cone angle $4 \pi$ and ending at cone angle $2 \pi$, so they depend
on the new work in \cite{HK3}.

Below we give a brief description of the boundary value problem involved
in this new version of the harmonic deformation theory, as well as some
applications to finite volume hyperbolic cone-manifolds.

We will assume, for simplicity, that $\M$ is a compact hyperbolic $3$-manifold
with a single torus as its boundary.  Hyperbolic manifolds with multiple torus
boundary components can be handled by using the same type of boundary
conditions on each one.  The theory extends to hyperbolic manifolds which
also have infinite volume geometrically finite ends, whose conformal structure
at infinity is assumed to be fixed, using the same techniques as in \cite{Br1}.

We further assume that the geometry near each boundary torus is 
modelled on the {\em complement} of an open tubular neighborhood of 
radius $R$ around
the singular set of a hyperbolic cone-manifold.  (A horospherical
neighborhood of a cusp is included by allowing $R= \infty$.)  In 
particular, we assume
that each torus has an intrinsic flat metric with constant principal 
normal curvatures
$\kappa$ and ${1 \over \kappa}$,  where $\kappa \geq 1$.
The normal curvatures and the tube radius, $R$, are related by
${\rm coth}~R = \kappa$ so they determine each other.  In fact, given
such a boundary structure, it can be canonically filled in.  In general,
the filled-in structure has ``Dehn surgery type singularities'' (see
\cite[Chap. 4]{thnotes}), which includes cone singularities with arbitrary cone
angle.  We say that $\M$ has {\em tubular boundary}.  This structure 
is described in
more detail in \cite{HK2}.

In \cite{HK1} specific closed $E$-valued $1$-forms, defined in a 
neighborhood of
the singular
locus, were exhibited which had the property that some complex linear 
combination of them
induced every possible infinitesimal change in the holonomy 
representation of a boundary
torus.  As a result, by standard cohomology theory, for any 
infinitesimal deformation
of the hyperbolic structure, it is possible to find a closed 
$E$-valued $1$-form $\hat \omega$
on $\M$ which equals such a complex linear combination of these 
standard forms in a neighborhood
of the torus boundary.  This combination of standard forms corresponds to the
terms $\omega_m + \omega_l$ in equation (\ref{omegadecomp}) in 
Section \ref{infharm}.

The standard forms are harmonic so the $E$-valued $1$-form $\hat \omega$
will be harmonic in a neighborhood of the boundary but not generally
harmonic on all of $\M$.  Since it represents a cohomology class in
$H^1(\M;E)$, it will be closed as an $E$-valued $1$-form, but it won't
generally be co-closed.  If we denote by $d_E$ and $\delta_E$
the exterior derivative and its adjoint on $E$-valued forms, then this means
that $d_E \hat \omega = 0$, but $\delta_E \hat \omega \neq 0$ in general.
(Note that $E$ is a flat bundle so that $d_E$ is the coboundary
operator for this cohomology theory.)

A representative for a cohomology class can be altered by a coboundary
without changing its cohomology class.  An $E$-valued $1$-form is a
coboundary precisely when it can be expressed as $d_E s$, where
$s$ is an $E$-valued $0$-form, i.e. a {\em global} section of $E$.
Thus, finding a harmonic (closed and co-closed) representative
cohomologous to $\hat \omega$ is equivalent to finding
a section $s$ such that
\begin{eqnarray}
\delta_E d_E s ~ = ~ - \delta_E \hat \omega.
\label{correction}
\end{eqnarray}
Then, $\omega = \hat \omega + d_E s$ satisfies
$\delta_E \omega = 0, \, d_E \omega = 0$; it is a closed and co-closed
(hence harmonic) representative in the same cohomology class as $\hat \omega$.

In \cite{HK1} it was shown that in order to solve equation (\ref{correction})
for $E$-valued sections, it suffices to solve it for the ``real
part" of $E$, where we are interpreting $E$ as the complexified
tangent bundle of $\M$ as discussed at the end of Section \ref{infharm}.
The real part of a section $s$ of $E$ is just a (real) section of the tangent
bundle of $\M$; i.e., a vector field, which we denote by $v$.  The
real part of $\delta_E d_E s$ equals $(\grad^* \grad + 2) \,v$, where
$\grad$ denotes the (Riemannian) covariant derivative and $\grad^*$
is its adjoint.  The composition $\grad^* \grad$ is sometimes called
the ``rough Laplacian" or the ``connection Laplacian".  We will denote
it by $\lap$.

To solve the equation $\delta_E d_E s = - \delta_E \hat \omega$,
we take the real part of $- \delta_E \hat \omega$, considered
as a vector field, and denote it by $\zeta$.  We then solve the equation
\begin{eqnarray} (\lap + 2) \, v ~=~ \zeta, \label{vfeqn}
\end{eqnarray}
for a vector field $v$ on $\M$.  As discussed in \cite{HK1}, this gives rise
to a section $s$ of $E$ whose real part equals $v$ which is a solution to
(\ref{correction}).  We denote the correction term $d_E s$ by
$\omega_c$.  In a neighborhood of the boundary $\hat \omega$
equals a combination of standard forms, $\omega_m + \omega_l$.
Thus, in a neighborhood of the boundary, we can decompose the harmonic
representative $\omega$ as $\omega = \omega_m + \omega_l + \omega_c$,
just as we did in (\ref{omegadecomp}) in Section \ref{infharm}.
Note that, since $\hat \omega$ was already harmonic in a neighborhood
of the boundary, the correction term $\omega_c$ will also be harmonic
on that neighborhood.

Because $\omega$ is harmonic, it will satisfy Weitzenb\"ock formulae
as described in Section \ref{effrigid}.  As will be outlined below we
will choose boundary conditions that will further guarantee that
$\omega = \eta + i *\! D\eta$ where $\eta$ is a $1$-form with values in
the tangent bundle of $\M$.  It decomposes as
$\eta = \eta_m + \eta_l + \eta_c$ in a neighborhood of the boundary and
satisfies the equation $D^*D\eta ~+~ \eta ~=~ 0$  on all of $\M$.
As before, taking an $L^2$ inner product and integrating by parts leads
to equation (\ref{bdryeq}).  The key to generalizing
the harmonic deformation theory from the previous sections is
finding boundary conditions that will guarantee that the contribution
to the boundary term of (\ref{bdryeq}) from the correction term
$\eta_c$ will be {\em non-positive}.  Once this is established, everything
else goes through without change.

In order to control the behavior of $d_E s = \omega_c$ (hence of $\eta_c$),
it is necessary to put restrictions on the domain
of the operator $(\lap + 2)$ in equation (\ref{vfeqn}).
On smooth vector fields with compact support the operator
$(\lap + 2)$ is self-adjoint.  It is natural to look for boundary conditions
for which self-adjointness still holds.  It is possible to find 
boundary conditions
which make this operator elliptic and self-adjoint with trivial kernel.
Standard theory then implies that equation (\ref{vfeqn}) is always 
uniquely solvable.

There are many choices for such boundary data.  Standard examples include the
conditions that either $v$ or its normal derivative be zero, 
analogous to Dirichlet
and Neumann conditions for the Laplacian on real-valued functions.
However, none of these standard choices of boundary data have the key property
that the Weitzenb\"ock boundary term for $\eta_c$ will necessarily be 
non-positive.

The main analytic result in \cite{HK3} is the construction of a boundary value
problem that has this key additional property.  We give a very brief 
description
of the boundary conditions involved.  (In particular, we will avoid discussion
of the precise function spaces involved.)

The first boundary condition on the vector fields allowed in the domain
of the operator $(\lap + 2)$ in equation (\ref{vfeqn}) is that its
($3$-dimensional) divergence vanish on the boundary.  This means that
the corresponding infinitesimal change of metric is volume preserving
at points on the boundary.  The combination of standard forms,
$\omega_m +\omega_l$, also induces infinitesimally volume preserving
deformations.  Once we prove the existence of a harmonic
$E$-valued $1$-form $\omega$ whose correction term comes from a vector
field satisfying this boundary condition, we can conclude that
$\omega$ is infinitesimally volume preserving at the boundary.
If we denote by ${\rm div}$ the function measuring the infinitesimal
change of volume at a point, then this means that ${\rm div}$ 
vanishes at the boundary.
However, for any harmonic $E$-valued $1$-form ${\rm div}$ satisfies the
equation:
\begin{eqnarray}
\lap~{\rm div} = -4 ~ {\rm div},
\end{eqnarray}
where $\lap$ here denotes the laplacian on functions given locally by the
sum of the {\em negatives} of the second derivatives.

A standard integration by parts argument shows that any function satisfying
such an equation and vanishing at the boundary must be identically zero.
Thus we can conclude that the deformation induced by $\omega$ is
infinitesimally volume preserving at {\em every} point in $\M$.  This allows
us to conclude that $\omega$ can be written as $\omega = \eta + i *\! D\eta$
where $\eta$ satisfies $D^*D\eta ~+~ \eta ~=~ 0$.  (See Proposition 2.6
in \cite{HK1}.)  The computation of the Weitzenb\"ock boundary term now
proceeds as before.

The second boundary condition is more complicated to describe.  Recall that
the boundary of $\M$ has the same structure as the boundary of a tubular
neighborhood of a (possibly singular) geodesic.  In particular, it has a
neighborhood which is foliated by tori which are equidistant from the 
boundary.
In a sufficiently small neighborhood, these surfaces are all embedded and,
on each of them, the nearest point projection to the boundary is a
diffeomorphism.  If we denote by $u$ the
($2$-dimensional) tangential component of the vector field $v$ at the
boundary, we can use these projection maps to pull back $u$ to these
equidistant surfaces.  We denote the resulting extension of $u$ to the
neighborhood of the boundary by $\hat u$.

In dimension $3$, the {\em curl} of a vector field is again a $3$-dimensional
vector field.  The second boundary condition is that the ($2$-dimensional)
tangential component of ${\rm curl}\,v$ agree with that of
${\rm curl}\,\hat u$ on the boundary.  Note that, on the boundary, the normal
component of ${\rm curl}\,v$ equals the curl of the $2$-dimensional vector
field $u$ (this curl is a function) so it automatically agrees with that
of $\hat u$.

As a partial motivation for this condition, consider a vector field
which generates an infinitesimal isometry in a neighborhood of the boundary
which preserves the boundary as a set.  Geometrically, it just translates
the boundary and all of the nearby equidistant surfaces along themselves.
Thus, this vector field is tangent to the boundary and to all the equidistant
surfaces.  It has the property that it equals the vector field $\hat u$ defined
above as the extension of its tangential boundary values.
Thus, the above condition can be viewed as an attempt to mirror properties
of infinitesimal isometries preserving the boundary.

To see why it might be natural to put conditions on the {\em curl} of
$v$, rather than on $v$ itself, consider the real-valued $1$-form $\tau$
dual (using the hyperbolic metric) to $v$.  Then, $\delta \tau$
and $*d\tau$ correspond, respectively, to the divergence and curl of $v$,
where $d$ denotes exterior derivative, $\delta$ its adjoint, and
$*$ is the Hodge star-operator.  Our boundary conditions can be viewed
as conditions on the exterior derivative and its adjoint applied
to this dual $1$-form.

However, the ultimate justification for these boundary conditions is
that they lead to a Weitzenb\"ock boundary term with the correct
properties, as long as the tube radius is sufficiently large.
A direct geometric proof of this fact is still lacking, as is an
understanding of the geometric significance of the value of the required
lower bound on the tube radius.  Nonetheless, the fact that the contribution
to the Weitzenb\"ock boundary from the correction term $\omega_c$
is always non-positive when it arises from a vector field satisfying
these boundary conditions can be derived by straightforward
(though somewhat intricate) calculation.  All the results from
the previous harmonic theory follow immediately.  For example,
we can conclude:

\begin{theorem}[\cite{HK3}] For a finite volume hyperbolic cone-manifold with
singularities along a link with tube radius at least
$\arctanh(1/\sqrt{3}) \approx 0.65848$,
there are no deformations of the hyperbolic structure keeping the 
cone angles fixed.
Furthermore, the nearby hyperbolic cone-manifold structures are
parametrized by their cone angles.
\label{bdrylocalparam}
\end{theorem}

\noindent In the statement of this theorem, when the singular link 
has more than one
component, having tube radius at least $R$ means that there are 
disjoint embedded tubes
of radius $R$ around all the components.

Besides being able to extend local rigidity rel cone angles and parametrization
by cone angles, the boundary value theory permits all of the estimates
involved in the effective rigidity arguments to go through.  In particular,
inequalities (\ref{eta0pos}), (\ref{normbound}), and (\ref{keybound}) 
continue to
hold.  The packing arguments require no restriction on cone angles so that
the proofs of the results on hyperbolic Dehn surgery (e.g., Theorems 
\ref{thm3} and
\ref{thm4}) go through unchanged.

In order to give an efficient description of the conclusions of these arguments
in this more general context we first extend some previous definitions.
Recall that, if $\M$ has a complete finite volume hyperbolic structure
with one cusp and $T$ is an embedded horospherical torus, the normalized
length of a simple closed curve $\gamma$ on $T$ is defined as the length
of the geodesic isotopic to $\gamma$ in the flat metric on $T$,
scaled to have area $1$.  A {\em weighted} simple closed curve is just
a pair, $(\lambda, \gamma)$, where $\lambda$ is a positive real number.
Its normalized length is then defined to be $\lambda$ times the normalized
length of $\gamma$.  If a basis is chosen for $H_1(T,\Z)$, the set of
isotopy classes of non-contractible simple closed curves on $T$ corresponds to
pairs $(p,q)$ of relatively prime integers.  Then a weighted simple
closed curve $(\lambda, \gamma)$ can be identified with the point
$(\lambda p, \lambda q) \in \R^2 \cong H_1(T,\Z) \otimes \R \cong
H_1(T,\R)$.  It is easy to check that the notion of normalized length
extends by continuity to any $(x,y) \in H_1(T,\R)$.

The {\em hyperbolic Dehn surgery space} for $\M$ (denoted $\DS (\M)$) is a
subset of $H_1(T,\R) \cup \infty$
which serves as a parameter space for (generally incomplete) 
hyperbolic structures
on $\M$ (with certain restrictions on the structure near its end).
In particular, if we view a weighted simple closed curve $(\lambda, \gamma)$
as an element of $H_1(T,\R)$, then saying that it is in $\DS (\M)$ means
that there is a hyperbolic cone-manifold structure on the manifold obtained
from $\M$ by doing Dehn filling with $\gamma$ as meridian which has the
core curve as the singular locus with cone angle $\frac {2 \pi} {\lambda}$.
The point at infinity corresponds to the complete hyperbolic structure on
$\M$.

Thurston's hyperbolic Dehn surgery theorem (see \cite{thnotes}) states
that $\DS (\M)$ always
contains an open neighborhood of $\infty$.  This, in particular, implies
that it contains all but a finite number of pairs $(p,q)$ of relatively prime
integers, which implies that all but a finite number of the manifolds obtained
by (topological) Dehn surgery are hyperbolic.  However, since most of these
pairs are clustered ``near" infinity, the statement that it contains an open
neighborhood of infinity is considerably stronger.  Again, Thurston's
proof is not effective; it provides no information about the size
of any region contained in hyperbolic Dehn surgery space.  The 
following theorem,
whose proof is analogous to that of Theorem \ref{thm3}
provides such information:

\begin{theorem}[\cite{HK3}] \label{shapethm}
Let $\M$ be a complete, finite volume orientable hyperbolic $3$-manifold with
one cusp, and let $T$ be a horospherical torus which is embedded as
a cross-section to the cusp of $\M$. Let
$(x,y)\in H_1(T,\R)$ have normalized length at least $7.583$.
Then there is a hyperbolic structure on $\M$ with Dehn surgery
coefficient $(x,y)$.  In particular, the hyperbolic Dehn surgery space
for $\M$ contains the complement of the ellipse around the origin determined
by the condition that the normalized length of $(x,y)$ is less than
$7.583$.  Furthermore, the volumes of hyperbolic structures in this
region differ from that of $\M$ by at most $0.306.$
\end{theorem}

\noindent {\bf Remark:} The homology group $H_1(T,\R)$ can be naturally
identified with the universal cover of $T$ so the flat metric on $T$,
normalized to have area $1$, induces a flat metric on $H_1(T,\R)$.
Then the ellipse in the above theorem becomes a metric disk of radius $7.583$.
  From this point of view, the theorem provides a universal size region
in $\DS (\M)$ (the complement of a ``round disk" of radius $7.583$), which
is {\em independent} of $\M$.  However, it is perhaps more interesting
to note that, if $H_1(T,\R)$ is more naturally identified with
$H_1(T,\Z) \otimes \R$, then this region actually reflects the
{\em shape} of $T$.

\medskip
Finally we note that these techniques also provide good estimates on 
the change in
geometry during hyperbolic Dehn filling as in Theorem \ref{shapethm}.
  For example, Schl\"afli's formula
(\ref{schlafli}) together with control on the length of the
singular locus (as in Theorem \ref{keyboundtheorem}) leads to
explicit upper and lower bounds for the
decrease in volume, $\Delta V$.
These bounds are independent of the cusped manifold $\M$, and
can be viewed as refinements of the asymptotic formula of
Neumann and Zagier \cite{neumann-zagier}:
$$\Delta V \sim {\pi^2 \over {L(x,y)}^2} \quad\text{~~as~~} L(x,y) 
\to \infty,$$
where $L(x,y)$ denotes the normalized length of the Dehn surgery coefficient
$(x,y) \in  H_1(T;\R)$. The details will appear in \cite{HK3}.

\end{document}